\documentclass[12pt,twoside]{amsart}
\usepackage{amssymb}

\usepackage[all]{xy}
\nonstopmode

\numberwithin{equation}{section} 

\font\tengothic=eufm10 scaled\magstep 1 \font\sevengothic=eufm7
scaled\magstep 1
\newfam\gothicfam
       \textfont\gothicfam=\tengothic
       \scriptfont\gothicfam=\sevengothic

\newtheorem{theorem}{Theorem}[section]

\newtheorem{proposition}[theorem]{Proposition}
\newtheorem{corollary}[theorem]{Corollary}
\newtheorem{conjecture}[theorem]{Conjecture}

\theoremstyle{definition}
\newtheorem{definition}[theorem]{Definition} 
\newtheorem{remark}[theorem]{Remark}

\def\mapright#1{\smash{ \mathop{\longrightarrow}
    \limits^{#1}}}

\newcommand{\cD}{{\mathcal D}}
\newcommand{\cA}{{\mathcal A}}

\newcommand{\cU}{{\mathcal U}}

\newcommand {\ZZ}{\mathbb{Z}}

\newcommand {\PP}{\mathbb{P}}

\begin{document}
\title[]{Betti numbers of determinantal ideals}

\author[ Rosa M.\ Mir\'o-Roig]{ Rosa M.\
Mir\'o-Roig$^{*}$}

\address{Facultat de Matem\`atiques,
Departament d'Algebra i Geometria, Gran Via de les Corts Catalanes
585, 08007 Barcelona, SPAIN } \email{miro@ub.edu}

\date{\today}
\thanks{$^*$ Partially supported by MTM2004-00666.}

\subjclass{Primary 13H15 13D02; Secondary 14M12}


\begin{abstract} Let $R=k[x_1, \cdots ,x_n]$ be a polynomial ring
and let $I\subset R$ be a graded ideal. In \cite{R}, R\"{o}mer
asked whether under the Cohen-Macaulay assumption the $i$-th Betti
number $\beta _{i}(R/I)$ can be bounded above by a function of the
maximal shifts in the minimal graded free $R$-resolution of $R/I$
as well as bounded below by a function of the minimal shifts. The
goal of this paper is to establish such bounds for graded
Cohen-Macaulay algebras  $k[x_1, \cdots ,x_n]/I$ when $I$ is a
standard determinantal ideal of arbitrary codimension. We also
discuss other examples as well as when these bounds are sharp.

\end{abstract}


\maketitle

\tableofcontents


  \section{Introduction} \label{intro}
Let $R=k[x_1, \cdots ,x_n]$ be a polynomial ring in $n$ variables
over a field $k$, let $deg(x_i)=1$ and let $I\subset R$ be a
graded ideal of arbitrary codimension. Consider the minimal graded
free $R$-resolution of $R/I$:

$$ 0 \longrightarrow  \oplus _{j\in \ZZ}R(-j)^{\beta_{p,j}(R/I)}\longrightarrow
\cdots  \longrightarrow   \oplus _{j\in
\ZZ}R(-j)^{\beta_{1,j}(R/I)} \longrightarrow R\longrightarrow R/I
\longrightarrow 0$$ where we denote $\beta_{i,j}(R/I)=\dim
Tor_{i}^R(R/I,k)_{j}$ the $(i,j)$-th graded Betti number of $R/I$
and $\beta _{i}(R/I)=\sum _{j\in \ZZ} \beta_{i,j}(R/I)$ is the
$i$-th total Betti number. Many important numerical invariants of
$I$ and the associated scheme can be read off from the minimal
graded free $R$-resolution of $R/I$. For instance, the Hilbert
polynomial, and hence the multiplicity $e(R/I)$ of $I$, can be
written down in terms of the shifts $j$ such that
$\beta_{i,j}(R/I)\ne 0$ for some $i$, $1\le i \le p$.

Let $c$ denote the codimension of $R/I$. Then $c\le p$ and
equality holds if and only if $R/I$ is Cohen-Macaulay. We define
$$m_i(I)=\text{min} \{ j\in \ZZ \mid \beta_{i,j}(R/I)\ne 0 \}$$ to
be the minimum degree shift at the $i$-th  step and
$$M_i(I)=\text{max} \{ j\in \ZZ \mid \beta_{i,j}(R/I)\ne 0 \}$$ to
be the maximum degree shift at the $i$-th  step. We will simply
write $m_i$ and $M_i$ when there is no confusion. If $R/I$ is
Cohen-Macaulay and has a pure resolution, i.e. $m_i=M_i$ for all
$i$, $1\le i \le c$, then Herzog and K\"{u}hl \cite{HK} and Huneke
and Miller \cite{HM} showed that
$$ e(R/I)=\frac{\prod_{i=1}^cm_i}{c!}$$
and $$ \beta _{i}(R/I)=(-1)^{i+1}\prod_{j\ne i}\frac{d_j}{d_j-d_i}
\text{ for } i=1,\cdots , c.$$ Since then there has been a
considerably effort to bound the multiplicity of a homogeneous
Cohen-Macaulay ideal $I\subset R$ in terms of the shifts in its
graded minimal $R$-free resolution; and  Herzog, Huneke and
Srinivasan have made the following conjecture ({\em minimal
conjecture}): If  $I\subset R$ be a graded Cohen-Macaulay ideal of
codimension $c$, then
$$\frac{\prod_{i=1}^cm_i}{c!}\le e(R/I)\le
\frac{\prod_{i=1}^cM_i}{c!}.$$

There is a growing body of the literature proving special cases of
the above conjecture. For example, it holds for complete
intersections \cite{HS}, powers of complete intersection ideals
\cite{GV}, perfect ideals with a pure resolution  \cite{HM},
perfect ideals with a quasi-pure resolution (i.e. $m_i\ge
M_{i-1}$) \cite{HS},  perfect ideals
 of codimension 2 \cite{HS}, Gorenstein ideals of codimension 3
 \cite{MNR} and
 standard determinantal ideals of arbitrary codimension \cite{M}.

Another natural question which naturally arises in this context is
whether under the Cohen-Macaulay assumption the $i$-th Betti
number $\beta _{i}(R/I)$ can be bounded above by a function of the
maximal shifts in the minimal graded free $R$-resolution of $R/I$
as well as bounded below by a function of the minimal shifts.
 In \cite{R}, R\"{o}mer made a natural guess
 \begin{equation}\label{guess1}
 \prod_{1\le j<i}\frac{m_j}{m_i-m_j}\prod_{i<j\le c}\frac{m_j}{m_j-m_i}
\le \beta_i(R/I)\le \prod_{1\le
j<i}\frac{M_j}{M_i-M_j}\prod_{i<j\le c}\frac{M_j}{M_j-M_i}
\end{equation} for $i=1,\cdots ,c$ and he showed that these bounds
hold if $I$ is a complete intersection or componentwise linear.
Moreover, in these cases we have equality above or below if and
only if $R/I$ has a pure resolution. Unfortunately, these bounds
are not always valid (see \cite{R}; Example 3.1). For
Cohen-Macaulay algebras with strictly quasi-pure resolution (i.e.,
$m_i>M_{i-1} $ for all $i$) he showed
\begin{equation}\label{guess2}
 \prod_{1\le j<i}\frac{m_j}{M_i-m_j}\prod_{i<j\le c}\frac{m_j}{M_j-m_i}
\le \beta_i(R/I)\le \prod_{1\le
j<i}\frac{M_j}{m_i-M_j}\prod_{i<j\le c}\frac{M_j}{m_j-M_i}
\end{equation} for $i=1,\cdots ,c$ and again we have equalities if and only if $R/I$ has
a pure resolution. Notice that $\prod_{1\le
j<i}\frac{M_j}{m_i-M_j}\prod_{i<j\le c}\frac{M_j}{m_j-M_i}$ may be
negative and thus, in general,  $\prod_{1\le
j<i}\frac{M_j}{m_i-M_j}\prod_{i<j\le c}\frac{M_j}{m_j-M_i}$  is
not a good candidate for being an upper bound. In \cite{R},
R\"{o}mer suggests as upper bound
\begin{equation}\label{guess3}
 \beta_i(R/I)\le \frac{1}{(i-1)!·(c-i)!}\prod_{j\ne i}M_j
\end{equation} for $i=1,\cdots ,c$ and he proved that the lower
bound in (\ref{guess2}) and the upper bound in (\ref{guess3}) hold
if $R/I$ is Cohen-Macaulay of codimension 2 and Gorenstein of
codimension 3.

\vskip 2mm It remains open if these bounds hold for other
interesting classes of ideals. The goal of this paper is to prove
that  the lower bound in (\ref{guess2}) and the upper bound in
(\ref{guess3}) work in the following classes of ideals
\begin{itemize}
\item standard determinantal ideals of arbitrary codimension $c$
(i.e.,
 ideals generated by the maximal minors of a $t\times (t+c-1)$ homogeneous
polynomial matrix), \item symmetric determinantal ideals defined
by the submaximal minors of a $t\times t$ homogeneous symmetric
matrix, \item determinantal ideals defined by the submaximal
minors of a $t\times t$ homogeneous matrix, and \item
arithmetically Cohen-Macaulay divisors on a variety of minimal
degree. \end{itemize}

Determinantal ideals are a central topic in both commutative
algebra and algebraic geometry. Due to their important role, their
study has attracted many researchers and has received considerable
attention in the literature.  Some of the most remarkable results
about determinantal ideals are due to  Eagon and  Hochster in
\cite{e-h}, and to  Eagon and  Northcott in \cite{e-n}. Eagon and
 Hochster proved that generic determinantal ideals are perfect.
 Eagon and  Northcott constructed a finite graded free
resolution for any determinantal ideal, and as a corollary, they
showed that determinantal ideals are perfect. Since then many
authors have made important contributions to the study of
determinantal ideals, and the reader can look at \cite{b-v},
\cite{BH}, \cite{KM} and \cite{eise} for background, history  and
a list of important papers.

\vskip 2mm \vskip 2mm Next we outline the structure of the paper.
In section 2, we first recall the basic facts on standard
determinantal ideals $I$
 of codimension $c$ defined by the maximal
minors of a $t\times (t+c-1)$ homogeneous matrix $\cA$ and the
associated complexes needed later on. We determine the minimal and
maximal shifts in the graded minimal free $R$-resolution of $R/I$
in terms of the degree matrix $\cU$ of $\cA$ and we prove that the
lower bound in (\ref{guess2}) and the upper bound in
(\ref{guess3}) work
 for standard determinantal ideals of arbitrary codimension $c$ (Theorem \ref{maintheorem}).
We also discuss when these bounds are sharp. After some
preliminaries, we devote section 3 to prove that the lower bound
in (\ref{guess2}) and the upper bound in (\ref{guess3}) work  for
determinantal (resp. symmetric determinantal) ideals defined by
the submaximal minors of a $t\times t$ homogeneous (resp.
symmetric) matrix (Theorems \ref{theorem2} and \ref{theorem3}). We
discuss other examples as well as when these bounds are sharp
(Theorem \ref{theorem4}).

 \vskip 4mm

\section{Standard Determinantal ideals}

In the first part of this  section, we  provide the background and
basic results on determinantal ideals needed in the sequel, and we
refer to \cite{b-v} and \cite{eise} for more details.

\vskip 2mm

Let $\cA$ be a homogeneous matrix, i.e. a matrix representing a
degree 0 morphism $\phi :F \longrightarrow G$ of free graded
$R$-modules. In this case, we denote by $I(\cA)$ the ideal of $R$
generated by the maximal minors of $\cA$ and by $I_j(\cA)$ the
ideal generated by the $j \times j$ minors of $\cA$.

\begin{definition}
 A homogeneous ideal $I\subset R$ of codimension $c$ is
called a \emph{standard determinantal} ideal if $I=I(\cA)$ for
some $t\times (t+c-1)$ homogeneous matrix $\cA$.
 \end{definition}

\vskip 2mm
 Let $I\subset R$  be a standard determinantal ideal
 of codimension    $c$
generated by the  maximal minors of a $t\times (t+c-1)$ matrix
$\cA=(f_{ij})_{i=1,...,t}^{j=1,...,t+c-1}$ where $f_{ij}\in {
k}[x_{1},...,x_{n}]$ are homogeneous polynomials of degree
$a_j-b_{i}$. The matrix $\cA$ defines a degree 0 map $$F=\oplus
_{i=1}^tR(b_i) \stackrel{\cA}{\longrightarrow}
G=\oplus_{j=1}^{t+c-1}R(a_j)$$ $$ v\mapsto v \cdot \cA$$ where
$v=(f_1,\cdots , f_t)\in F$ and we assume without loss of
generality that $\cA$ is minimal; i.e., $f_{ij}=0$ for all $i,j$
with $b_{i}=a_{j}$. If we let $u_{i,j}=a_j-b_i$ for all $j=1,
\dots , t+c-1$ and $i=1, \dots , t$, the matrix
$\cU=(u_{i,j})_{i=1,...,t}^{j=1,...,t+c-1}$ is called the {\em
degree matrix} associated to $I$. By re-ordering degrees, if
necessary,   we may also assume that $b_1 \ge ... \ge b_t$ and
$a_1 \le a_2\le ... \le a_{t+c-1}$. In particular, we have:

\begin{equation} \label{order} u_{i,j}\le u_{i+1,j} \quad \text{ and }
\quad  u_{i,j}\le u_{i,j+1} \text{ for all } i,j.
\end{equation}

 \vskip 2mm  Note that the
degree matrix $\cU$ is completely determined by  $u_{1,1}$,
$u_{1,2}$, ... , $u_{1,c}$,   $u_{2,2}$, $u_{2,3}$, ... ,
$u_{2,c+1}$, ..., $u_{t,t}$, $u_{t,t+1}$, ... , $u_{t,c+t-1}$
because of the identity
$u_{i,j}+u_{i+1,j+1}-u_{i,j+1}-u_{i+1,j}=0$ for all $i,j$.
Moreover, the graded Betti numbers in the minimal free
$R$-resolution of
 $R/I(\cA)$ depend only upon the integers
 $$\{ u_{i,j} \} _{1\le i \le t}^{i\le j \le c+i-1}\subset \{u_{i,j}\}_{i=1,...t}^{j=1,...,t+c-1}$$ as described
 below.

\vskip 4mm
\begin{proposition}\label{miMi}
 Let $I\subset R$ be a standard determinantal ideal of codimension $c$ with
 degree matrix $\cU=(u_{i,j})_{i=1,...,t}^{j=1,...,t+c-1}$ as above.
 Then we have:
 \vskip 2mm
 \begin{itemize}
 \item[(1)] $m_i=u_{1,1}+u_{1,2}+\cdots +u_{1,i}+u_{2,i+1}+u_{3,i+2}+\cdots +u_{t,t+i-1} $
  for  $1\le i \le c$,

  \vskip 2mm
 \item[(2)] $M_i=  u_{1,c-i+1}+u_{2,c-i+2}+\cdots +u_{t,t+c-i}+u_{t,t+c-i+1}+u_{t,t+c-i+2}
 +\cdots +u_{t,t+c-1} $ for  $1\le i \le c$,

 \vskip 2mm
 \item[(3)] $\beta _{i}(R/I)= {t+c-1\choose t+i-1}{t+i-2\choose i-1}$
 for $i=1,\cdots , c$.
 \end{itemize}
 \end{proposition}

\begin{proof}  We denote by
$\varphi :F \longrightarrow G$ the morphism of free graded
$R$-modules of rank $t$ and $t+c-1$, defined by the homogeneous
matrix $\cA$ associated to $I$. The Eagon-Northcott complex
${\cD}_0(\varphi^*):$

$$0 \longrightarrow \wedge^{t+c-1}G^* \otimes S_{c-1}(F)\otimes \wedge^tF\longrightarrow
\wedge^{t+c-2} G ^*\otimes S _{c-2}(F)\otimes \wedge
^tF\longrightarrow \ldots \longrightarrow$$ $$\wedge^{t}G^*
\otimes S_{0}(F)\otimes \wedge^tF\longrightarrow R \longrightarrow
R/I \longrightarrow 0 $$

\vskip 2mm \noindent gives us a graded minimal free $R$-resolution
of $R/I$ (See, for instance \cite{b-v}; Theorem 2.20 and
\cite{eise}; Corollary A2.12 and Corollary A2.13). Now the result
follows after a straightforward computation taking into account
that  $$\wedge^tF =R(\sum_{i=1}^t b_{i}),$$ $$  S_a(F)=\bigoplus
_{1\le i_1\le \cdots \le i_a\le t}R(\sum_{j=1}^{a}b_{i_j}),$$
$$G^*=\bigoplus _{j=1}^{t+c-1}R(-a_{j}) \quad \text{and} $$ $$
\wedge^b G^*=\bigoplus _{1\le i_1<\cdots <i_b\le
t+c-1}R(-\sum_{j=1}^{b}a_{i_{j}}).$$
\end{proof}

\begin{remark} \label{keyremark} Let $I\subset R$ be a standard determinantal ideal.
It is worthwhile to point out that the $i$-th total Betti number
$\beta _i(R/I)$ in the minimal free $R$-resolution of $R/I$ depend
only upon the size $t\times (t+c-1)$ of the homogeneous matrix
$\cA$ associated to $I$.
\end{remark}

We are now ready to state the main result of this short note.

\vskip 2mm
\begin{theorem}\label{maintheorem} Let $I\subset R$ be a
standard determinantal ideal of codimension $c$. Then, we have:
\begin{equation}
 \prod_{1\le j<i}\frac{m_j}{M_i-m_j}\prod_{i<j\le c}\frac{m_j}{M_j-m_i}
\le \beta_i(R/I)\le \frac{1}{(i-1)!·(c-i)!}\prod_{j\ne i}M_j
\end{equation} for $i=1,\cdots ,c$.
In addition, the bounds are reached for all $i$ if and only if
$R/I$ has a pure resolution if and only if $u_{i,j}=u_{r,s}$ for
all $1\le i,r\le t$ and $1\le j,s\le t+c-1$.
\end{theorem}
\begin{proof} We will first prove the result for $J$ being $J\subset R$  a
standard determinantal ideal of codimension $c$ whose associated
matrix $\cA$ is a $t\times (t+c-1)$ matrix with all its entries
linear forms. In this case, for all $i=1,\cdots ,c$, we have (see
Proposition \ref{miMi}): $$m_i(J)=M_i(J)=t+i-1 \text{ and } \beta
_i(R/J)= {t+c-1\choose t+i-1}{t+i-2\choose i-1}.$$
  Therefore, $R/J$ has a pure
 resolution and it follows from \cite{HK} and \cite{HM} that
 \begin{eqnarray*} \prod_{1\le j<i}\frac{m_j(J)}{M_i(J)-m_j(J)}\prod_{i<j\le
 c}\frac{m_j(J)}{M_j(J)-m_i(J)}& = &  \prod_{1\le j<i}\frac{t+j-1}{i-j}\prod_{i<j\le
 c}\frac{t+j-1}{j-i}\\ & = & \beta _i(R/J) \\ & = &
\prod_{1\le j<i}\frac{t+j-1}{i-j}\prod_{i<j\le
 c}\frac{t+j-1}{j-i}\\
 & = & \prod_{1\le
j<i}\frac{M_j(J)}{m_i(J)-M_j(J)}\prod_{i<j\le c}\frac{M_j(J)}{m_j(J)-M_i(J)}\\
& = & \frac{1}{(i-1)!·(c-i)!}\prod_{j\ne i}M_j(J).
 \end{eqnarray*}

We will now prove the general case. Let $I$ be a standard
determinantal ideal of codimension $c$ with associated degree
matrix $\cU=(u_{i,j})_{i=1,...,t}^{j=1,...,t+c-1}$. Since, for all
$i=1, \cdots , c$, we have $$M_i(I)\ge m_i (I) \ge
t+i-1=m_i(J)=M_i(J),$$ it follows from Proposition \ref{miMi} (3)
and Remark \ref{keyremark}  that  \begin{eqnarray*}\beta _i(R/I) &
= & \beta _i(R/J)
\\ & = & \prod_{1\le j<i}\frac{t+j-1}{i-j}\prod_{i<j\le
 c}\frac{t+j-1}{j-i}\\ & = & \frac{1}{(i-1)!·(c-i)!}\prod_{j\ne
i}M_j(J)
 \\ & \le & \frac{1}{(i-1)!·(c-i)!}\prod_{j\ne
i}M_j(I) \end{eqnarray*} and this completes the proof of the upper
bound.

Let us now prove the lower bound. Recall that if $r\ge m$ and
$s\ge n$ then $u_{r,s} \ge u_{m,n}$. Therefore, for $1\le j<i$, we
get, using Proposition \ref{miMi}, that
\begin{eqnarray*} m_j & = & u_{1,1}+u_{1,2}+\cdots +u_{1,j}+u_{2,j+1}+u_{3,j+2}+\cdots
+u_{t,t+j-1} \\ & \le & (t+j-1)u_{t,t+j-1},
\end{eqnarray*}
\begin{eqnarray*}
M_i-m_j & = &  u_{1,c-i+1}+u_{2,c-i+2}+\cdots
+u_{t,t+c-i}+u_{t,t+c-i+1}+u_{t,t+c-i+2}
 +\cdots +u_{t,t+c-1} \\
 & - & (u_{1,1}+u_{1,2}+\cdots +u_{1,j}+u_{2,j+1}+u_{3,j+2}+\cdots
+u_{t,t+j-1})\\
& \ge & u_{t,t+c-i+j}+u_{t,t+c-i+j+1}+\cdots +u_{t,t+c-1}\\ & \ge
& (i-j)u_{t,t+c-i+j},
\end{eqnarray*} and
\begin{eqnarray*} \frac{m_j}{M_i-m_j}  & \le &
\frac{(t+j-1)u_{t,t+j-1}}{(i-j)u_{t,t+c-i+j}} \\
& \le & \frac{t+j-1}{i-j}.
\end{eqnarray*}

For $i<j\le c$, we have
\begin{eqnarray*} m_j & = & u_{1,1}+u_{1,2}+\cdots +u_{1,j}+u_{2,j+1}+u_{3,j+2}+\cdots
+u_{t,t+j-1} \\
& = & \left\{
\begin{array}{ll}
    u_{1,1}+u_{1,2}+\cdots +u_{1,t+i-1}+\cdots +u_{1,j}+u_{2,j+1}+\cdots
+u_{t,t+j-1}, & \mbox{ if } t+i\le j, \\
   u_{1,1}+\cdots +u_{1,j}+u_{2,j+1}+\cdots
   +u_{t+i-j,t+i-1}+\cdots
+u_{t,t+j-1} , & \mbox{ if } t+i> j;
\end{array}
\right.
\end{eqnarray*}
\begin{eqnarray*}
M_j-m_i & = &  u_{1,c-j+1}+u_{2,c-j+2}+\cdots
+u_{t,t+c-j}+u_{t,t+c-j+1}+u_{t,t+c-j+2}
 +\cdots +u_{t,t+c-1} \\
 & - & (u_{1,1}+u_{1,2}+\cdots +u_{1,i}+u_{2,i+1}+u_{3,i+2}+\cdots
+u_{t,t+i-1})\\
& \ge & u_{t,t+c+i-j}+u_{t,t+c+i-j+1}+\cdots +u_{t,t+c-1} \\
& \ge & \left\{
\begin{array}{ll}
    u_{1,t+i}+\cdots +u_{1,j}+u_{2,j+1}+\cdots
+u_{t,t+j-1}, & \mbox{ if } t+i\le j, \\
   u_{t+i-j+1,t+i}+\cdots
+u_{t,t+j-1} , & \mbox{ if } t+i> j.
\end{array}
\right.
\end{eqnarray*}

Therefore, if $i<j\le c$ and $t+i\le j$, we get
\begin{eqnarray*} \frac{m_j}{M_j-m_i}  & \le & \frac{u_{1,1}+u_{1,2}+\cdots +u_{1,t+i-1}}{
u_{1,t+i}+\cdots +u_{1,j}+u_{2,j+1}+\cdots
+u_{t,t+j-1}}+1 \\
& \le & \frac{(t+i-1)u_{1,t+i-1}}{(j-i)u_{1,t+i}}+1\\
& \le & \frac{t+i-1}{j-i}+1 \\ & = & \frac{t+j-1}{j-i},
\end{eqnarray*}
and if $i<j\le c$ and $t+i> j$, we get
\begin{eqnarray*} \frac{m_j}{M_j-m_i}  & \le & \frac{u_{1,1}+\cdots +u_{1,j}+u_{2,j+1}+\cdots +
u_{t+i-j,t+i-1}}{ u_{t+i-j+1,t+i}+\cdots +u_{t,t+j-1}
}+1\\
& \le & \frac{(t+i-1)u_{t+i-j,t+i-1}}{(j-i)u_{t+i-j+1,t+i}}+1\\
& \le & \frac{t+i-1}{j-i}+1 \\ & = & \frac{t+j-1}{j-i}.
\end{eqnarray*}

 Hence, putting altogether, we obtain
\begin{eqnarray*}
\prod_{1\le j<i}\frac{m_j}{M_i-m_j}\prod_{i<j\le
c}\frac{m_j}{M_j-m_i} & \le  & \prod_{1\le
j<i}\frac{t+j-1}{i-j}\prod_{i<j\le
 c}\frac{t+j-1}{j-i}\\ & =  & \beta _i(R/J) \\ & = & \beta _i(R/I)
\end{eqnarray*} and this completes the proof of the lower bound.
 Checking the inequalities we easily see that we have equality
above and below for all $1\le i \le c$ if and only if $R/I$ has a
pure resolution. This concludes the proof of the Theorem.
\end{proof}

\begin{remark} Since a complete intersection ideal $I$ of
arbitrary codimension and Cohen-Macaulay ideals of codimension 2
are examples of standard determinantal ideals, we recover
\cite{R}; Theorem 2.1 and Corollary 4.2.
\end{remark}

Since the power $I^s$ of a complete intersection ideal $I\subset
R$ is an example of standard determinantal ideal,  as a corollary
of Theorem \ref{maintheorem}, we have

\begin{corollary}\label{coro1} Let $I\subset R$ be a complete
intersection ideal of codimension $c$ and let $s$ be any positive
integer. Then, it holds
\begin{equation}
 \prod_{1\le j<i}\frac{m_j}{M_i-m_j}\prod_{i<j\le c}\frac{m_j}{M_j-m_i}
\le \beta_i(R/I^s)\le \frac{1}{(i-1)!·(c-i)!}\prod_{j\ne i}M_j
\end{equation} for $i=1,\cdots ,c$.
\end{corollary}

\section{Ideals defined by submaximal minors}

\vskip 2mm
 The first goal of this section is to prove that
   the lower bound in (\ref{guess2}) and the upper bound in
(\ref{guess3}) work  for k-algebras $k[x_1,\cdots ,x_n]/I$ being
$I$ a perfect ideal generated by the submaximal minors of a $t
\times t$ homogeneous symmetric matrix. A classical homogeneous
ideal that can be generated by the submaximal minors of a $t
\times t$ homogeneous symmetric matrix is the ideal of the
Veronese surface $X\subset \PP^5$. Indeed, the ideal of the
Veronese surface $X\subset \PP^5=Proj(k[x_1,\cdots ,x_6])$ can be
generated by the $2\times 2$ minors of the symmetric matrix
$$\begin{pmatrix}
x_6 & x_1 & x_2 \\
x_1 & x_3  & x_4 \\
x_2 & x_4 & x_5
\end{pmatrix}.
$$

Let us now fix some notation. Let $I\subset S=k[x_1, \cdots ,x_n]$
 be a codimension 3, perfect ideal generated by the submaximal
minors of a $t \times t$ homogeneous symmetric matrix
$\cA=(f_{ji})_{i,j=1,...,t}$ where $f_{ji}\in {
k}[x_{1},...,x_{n}]$ are homogeneous polynomials of degree
$a_i+a_{j}$, i.e. $I=I_{t-1}(\cA)$. We denote by

$$\cU=\begin{pmatrix}
2a_1 & a_1+a_2 & \cdots & a_1+a_t \\
a_1+a_2 & 2a_2 & \cdots & a_2+a_t \\
\vdots &  \vdots & &  \vdots\\
 a_1+a_t &
a_2+a_t & \cdots & 2a_t
\end{pmatrix}
$$

\vskip 2mm \noindent the degree matrix of $\cA $. We assume that
$a_1\le a_2 \le \cdots \le a_t$. The determinant of $\cA $ is a
homogeneous polynomial of degree $\ell :=2(a_1+a_2+\cdots +a_t)$.
Note that $a_i+a_j$ is an integer for all $1\le i \le j \le t$
while $a_i$ does not necessarily need to be an integer.

 \vskip 2mm  Note that the
degree matrix $\cU$ is completely determined by  $a_{1}$,  ... ,
$a_{t}$. Moreover, the graded Betti numbers in the minimal free
$S$-resolution of
 $S/I_{t-1}(\cA)$ depend only upon the integers  $a_{1}$,  ... ,
 $a_t$
as we will describe now. To this end, we recall Jozefiak's result
about the resolution of ideals generated by minors of a symmetric
matrix.

\vskip 2mm Let $R$ be a commutative ring with identity and let
$X=(x_{ij})$ be a symmetric $t \times t$ matrix with entries in
$R$. Write $Y=(y_{ij})$ for the matrix of cofactors of $X$, i.e.,
$y_{ij}=(-1)^{i+1}X^{i}_j$ where $X^{i}_j$ stands for the minor of
$X$ obtained by deleting the $i$-th column and the $j$-th row of
$X$. The matrix $Y$ is also a symmetric matrix. Let $M_t(R)$ be
the free $R$-module of all $t\times t$ matrices over $R$ and let
$A_t(R)$ be the free $R$-submodule of $M_t(R)$ consisting of all
alternating matrices. Denote by $tr:M_t(R)\longrightarrow R$ the
trace map. By \cite{J}; Theorem 3.1, the free complex of length 3
associated to $X$:

$$0\longrightarrow A_t(R)\mapright{d_3} Ker(M_t(R)\mapright{tr}   R)
\mapright{d_2}
   M_t(R)/A_t(R) \mapright{d_1}  R$$

\noindent where the corresponding differentials are defined as
follows:

\vskip 2mm $d_3(A)=AX$,

$d_2(N)=XN$ mod $A_t(R)$, and

$d_1(M$mod $A_t(R))=tr(YM)$

\vskip 2mm \noindent is acyclic and gives a free resolution of
$R/I_{t-1}(X)$. So, we obtain

\vskip 4mm
\begin{proposition}\label{miMi2}
 Let $I\subset S=k[x_1, \cdots ,x_n]$ be  a perfect
ideal of codimension $3$ generated by the submaximal minors of a
symmetric matrix $\cA$. Let  $$\cU=\begin{pmatrix}
2a_1 & a_1+a_2 & \cdots & a_1+a_t \\
a_1+a_2 & 2a_2 & \cdots & a_2+a_t \\
\vdots &  \vdots & &  \vdots\\
 a_1+a_t &
a_2+a_t & \cdots & 2a_t
\end{pmatrix}
$$ be the degree matrix and $\ell :=2(a_1+a_2+\cdots +a_t)$.
 Then, we have:
 \vskip 2mm
\begin{itemize}
 \item[(1)] $m_1=\ell-2a_t$ and $M_1=\ell-2a_1$,
\vskip 2mm  \item[(2)]
 $m_2=\ell-a_t+a_1$ and $M_2=\ell-a_1+a_t$,
 \vskip 2mm \item[(3)]
 $m_3=\ell+a_1+a_2$ and $M_3=\ell+a_{t-1}+a_t$, and
 \vskip 2mm \item[(4)] $\beta_1(R/I)= {t+1\choose 2}$,  $\beta_2(R/I)=t^2-1 $, and $\beta_3(R/I)={t\choose 2}
 $.
 \end{itemize}
 \end{proposition}

\begin{proof}
By \cite{J}; Theorem 3.1, $I$ has a minimal free $S$-resolution of
the following type:

\begin{equation} \label{jozefiak} 0\longrightarrow
\oplus _{1\le i <j\le t} S(-a_i-a_j-\ell)\longrightarrow (\oplus
_{1\le i,j\le t}S(-\ell-a_i+a_j))/S(-\ell) \longrightarrow
\end{equation} $$ \oplus _{1\le i\le j\le t}S(a_i+a_t-\ell
)\longrightarrow I\longrightarrow 0.$$ So,  the maximum and
minimum degree shifts at the $i$-th  step are
 \begin{itemize}
 \item[(1)] $m_1=\ell-2a_t$ and $M_1=\ell-2a_1$,
 \item[(2)] $m_2=\ell-a_t+a_1$ and $M_2=\ell-a_1+a_t$,
  \item[(3)] $m_3=\ell+a_1+a_2$ and $M_3=\ell+a_{t-1}+a_t$.
 \end{itemize} and the $i$-the total Betti numbers are:
\begin{itemize}
 \item[(4)]  $\beta_1(R/I)={t+1\choose 2} $,  $\beta_2(R/I)=t^2-1 $, and $\beta_3(R/I)={t\choose 2} $
 \end{itemize}
 which proves what we want.
\end{proof}

\begin{remark} \label{keyremark2} Let $I\subset R$ be a perfect ideal of codimension $3$ generated
by the submaximal minors of a symmetric matrix. It is worthwhile
to point out that the $i$-th total Betti number $\beta _i(R/I)$ in
the minimal free $R$-resolution of $R/I$ depend only upon the size
$t\times t$ of the homogeneous symmetric matrix $\cA$ associated
to $I$.
\end{remark}

We are now ready to bound the $i$-the total Betti number of
 codimension 3, perfect ideals generated
by the submaximal minors of a symmetric matrix in terms of the
shifts in its minimal free $R$-resolution.

\begin{theorem} \label{theorem2}
Let $I\subset S$ be a perfect ideal of codimension $3$ generated
by the submaximal minors of a $t\times t$ symmetric matrix.
 Then, we have:
\begin{equation}
 \prod_{1\le j<i}\frac{m_j}{M_i-m_j}\prod_{i<j\le 3}\frac{m_j}{M_j-m_i}
\le \beta_i(R/I)\le \frac{1}{(i-1)!·(c-i)!}\prod_{j\ne i}M_j
\end{equation} for $1\le i\le 3$.
In addition, the bounds are reached for all $i$ if and only if
$R/I$ has a pure resolution if and only if $a_1=\cdots = a_t$.
\end{theorem}
\begin{proof} We will first prove the result for $J$ being $J\subset R$  a
codimension $3$ perfect ideal generated by the submaximal minors
of a $t\times t$ symmetric matrix $\cA$  with  linear entries. In
this case, for all $1\le i\le 3$, we have (see Proposition
\ref{miMi2}) $$m_i(J)=M_i(J)=t+i-2,$$ $$\beta_1(R/J)= {t+1\choose
2}, \quad \beta_2(R/J)=t^2-1 \text{ and } \beta_3(R/J)={t\choose
2}.$$
 Therefore, $R/J$ has a pure
 resolution and it follows from \cite{HK} and \cite{HM} that
 \begin{eqnarray*} \prod_{1\le j<i}\frac{m_j(J)}{M_i(J)-m_j(J)}\prod_{i<j\le
 3}\frac{m_j(J)}{M_j(J)-m_i(J)}& = &  \prod_{1\le j<i}\frac{t+j-2}{i-j}\prod_{i<j\le
 3}\frac{t+j-2}{j-i}\\ & = & \beta _i(R/J) \\ & = &
\prod_{1\le j<i}\frac{t+j-2}{i-j}\prod_{i<j\le
 3}\frac{t+j-2}{j-i}\\
 & = & \prod_{1\le
j<i}\frac{M_j(J)}{m_i(J)-M_j(J)}\prod_{i<j\le 3}\frac{M_j(J)}{m_j(J)-M_i(J)}\\
& = & \frac{1}{(i-1)!·(3-i)!}\prod_{j\ne i}M_j(J).
 \end{eqnarray*}

We will now prove the general case. Let $I$ be a perfect ideal of
codimension $3$ generated by the submaximal minors of a $t\times
t$ symmetric matrix, let  $$\cU=\begin{pmatrix}
2a_1 & a_1+a_2 & \cdots & a_1+a_t \\
a_1+a_2 & 2a_2 & \cdots & a_2+a_t \\
\vdots &  \vdots & &  \vdots\\
 a_1+a_t &
a_2+a_t & \cdots & 2a_t
\end{pmatrix}
$$ be its degree matrix and $\ell :=2(a_1+a_2+\cdots +a_t)$.
 Since, for all $1\le i \le 3$, we have
$$M_i(I)\ge m_i (I) \ge t+i-2=m_i(J)=M_i(J),$$ it follows from
Proposition \ref{miMi2} (3) and Remark \ref{keyremark2} that
\begin{eqnarray*}\beta _i(R/I) & = & \beta _i(R/J)
\\ & = & \prod_{1\le j<i}\frac{t+j-2}{i-j}\prod_{i<j\le
 3}\frac{t+j-2}{j-i}\\ & = & \frac{1}{(i-1)!·(3-i)!}\prod_{j\ne
i}M_j(J)
 \\ & \le & \frac{1}{(i-1)!·(3-i)!}\prod_{j\ne
i}M_j(I) \end{eqnarray*} and this completes the proof of the upper
bound.

Let us now prove the lower bound. Using again Proposition
\ref{miMi2} and Remark \ref{keyremark2}, we have
\begin{eqnarray*} \frac{m_1}{M_3-m_1}·\frac{m_2}{M_3-m_2} & = &
\frac{\ell
-2a_t}{3a_t+a_{t-1}}·\frac{\ell-a_t+a_1}{a_{t-1}+2a_t-a_1} \\
& \le & \frac{\ell
-2a_t}{4a_{t-1}}·\frac{\ell-a_t+a_1}{2a_t} \\
& \le & \frac{2(t-1)a_{t-1}}{4a_{t-1}}·\frac{2ta_t}{2a_t} \\
& = & \frac{(t-1)t}{2}\\
& = & \beta _3(R/J) \\
& = & \beta _3(R/I).
\end{eqnarray*}
Analogously, we obtain
$$\frac{m_1}{M_2-m_1}·\frac{m_3}{M_3-m_2}  \le (t-1)(t+1) = \beta _2 (R/J)= \beta
_2 (R/I)$$ and $$\frac{m_2}{M_2-m_1}·\frac{m_3}{M_3-m_1} \le
\frac{t(t+1)}{2} = \beta _1 (R/J)= \beta _1(R/I)$$
 and this completes the proof of the lower
bound. It is easy to see that the inequality is an equality for
all $i$ if and only if $a_1=a_2\cdots =a_t$ if and only if $R/I$
has a pure resolution.
\end{proof}

Now, we will focus our attention on
 determinantal
 ideals $I$ generated by the submaximal minors
 of a $t\times t$ quadratic homogeneous matrix. By  \cite{GN}; Th\'{e}or\`{e}me 2,
 $I$ is a Gorenstein ideal of codimension 4 and we will prove
  the lower bound in (\ref{guess2}) and the upper bound in
(\ref{guess3})  for
 such kind of perfect ideals.

 To this end, let
$\cA=(f_{ji})_{i,j=1,...,t}$ be a homogeneous quadratic matrix
with entries    homogeneous polynomials $f_{ji}\in {
k}[x_{1},...,x_{n}]$ of degree $a_j-b_{i}$. We assume without loss
of generality that $\cA$ is minimal; i.e., $f_{ji}=0$ for all
$i,j$ with $b_{i}=a_{j}$. If we let $u_{j,i}=a_j-b_i$ for all
$1\le i,j \le t$, the matrix $\cU=(u_{ji})_{i,j=1,...t}$ is called
the {\em degree matrix} associated to $I$.
 By  re-ordering
degrees, if necessary,   we may also assume that $b_1 \le ... \le
b_t$ and $ a_1\le ... \le a_{t}$. The determinant of $\cA $ is a
homogeneous polynomial of degree
$$s: =\deg(\det(\cA))=\sum_{j=1}^ta_{j}-\sum_{i=1}^tb_i.$$

 \vskip 2mm   Let
us now compute the $i$-th total Betti numbers in the minimal free
$R$-resolution of
 $R/I_{t-1}(\cA)$ in terms of the degree matrix
 $\{u_{j,i}\}_{i,j=1,...,t}$ of $\cA$.

\vskip 4mm
\begin{proposition}\label{miMi3}
 Let $I\subset R$ be a determinantal ideal of codimension $4$ generated by the
 submaximal minors of $t\times t$ quadratic homogeneous matrix
 $\cA$ with
 degree matrix $\cU=(u_{j,i})_{i,j=1,...,t}$ with
 $u_{j,i}=a_j-b_i$. Set $s: =\sum_{j=1}^ta_{j}-\sum_{i=1}^tb_i.$
 Then we have:
 \vskip 2mm
 \begin{itemize}
 \item[(1)]  $m_1(I)=s+b_1-a_t$, $M_1(I)=s+b_t-a_1$,
 \item[(2)] $m_2(I)=min(s+b_1-b_t,s+a_1-a_t)$, $M_2(I)=max(s+b_t-b_1,s+a_t-a_1)$,
 \item[(3)]  $m_3(I)=s-b_t+a_1$, $M_3(I)=s+a_t-b_1$,
 \item[(4)] $m_4(I)=M_4(I)=2s$, and
 \item[(5)]  $\beta_1(R/I)=t^2 $,  $\beta_2(R/I)=t^2-2 $,
 $\beta_3(R/I)=t^2
 $ and $\beta _4(R/I)=1$.
 \end{itemize}
 \end{proposition}
 \begin{proof}
 We denote by $$F:=\oplus _{i=1}^tR(b_i)\stackrel
{\cA}{ \longrightarrow} G:=\oplus _{j=1}^tR(a_j)$$ the morphism of
free graded $R$-modules of rank $t$, defined by the homogeneous
matrix $\cA$ associated to $I$. By \cite{GN}; Th\'{e}or\`{e}me 2,
$R/I$ has a minimal graded free $R$-resolution of the following
type:

\begin{equation}\label{gulliksen}0 \longrightarrow R(-2s)
\longrightarrow \bigoplus _{1\le i,j \le j}R(b_j-a_i-s)
\longrightarrow \end{equation} $$(\bigoplus _{1\le i,j\le t \atop
i\ne j}R(a_i-a_j-s) )\oplus ( \bigoplus _{1\le i,j\le t \atop i\ne
j}R(b_i-b_j-s))\oplus R(-s)^{2t-2}\longrightarrow \bigoplus _{1\le
i,j \le j}R(a_j-b_i-s) \longrightarrow I  \longrightarrow 0
$$
 and, hence,  the maximum and
minimum degree shifts at the $i$-th  step are
 \begin{itemize}
 \item[(1)]  $m_1=s+b_1-a_t$, $M_1=s+b_t-a_1$,
 \item[(2)] $m_2=min(s+b_1-b_t,s+a_1-a_t)$, $M_1=max(s+b_t-b_1,s+a_t-a_1)$,
 \item[(3)]  $m_3=s-b_t+a_1$, $M_3=s+a_t-b_1$, and
 \item[(4)] $m_4=M_4=2s$
 \end{itemize} and the $i$-the total Betti numbers are:
\begin{itemize}
  \item[(5)]  $\beta_1(R/I)=t^2 $,  $\beta_2(R/I)= t^2-2$,
  $\beta_3(R/I)=t^2
 $ and $\beta _4(R/I)=1$
 \end{itemize}
 which proves what we want.
\end{proof}

\begin{remark} \label{keyremark3} Let $I\subset R$ be a perfect ideal of codimension $4$ generated
by the submaximal minors of a homogeneous $t\times t$ matrix. It
is worthwhile to point out that the $i$-th total Betti number
$\beta _i(R/I)$ in the minimal free $R$-resolution of $R/I$ depend
only on $t$.
\end{remark}

Arguing as in Theorems \ref{maintheorem} and \ref{theorem2} we can
bound the $i$-th total Betti number of a
 codimension 4, Gorenstein ideal generated
by the submaximal minors of a quadratic matrix in terms of the
shifts in its minimal free $R$-resolution and we get:

\begin{theorem} \label{theorem3} Let $I\subset R$ be a determinantal ideal of codimension $4$ generated by the
 submaximal minors of $t\times t$ quadratic homogeneous matrix
 $\cA$.
 Then, we have:
\begin{equation}
 \prod_{1\le j<i}\frac{m_j}{M_i-m_j}\prod_{i<j\le 4}\frac{m_j}{M_j-m_i}
\le \beta_i(R/I)\le \frac{1}{(i-1)!·(4-i)!}\prod_{j\ne i}M_j
\end{equation} for $1\le i\le 4$.
In addition, the bounds are reached for all $i$ if and only if
$R/I$ has a pure resolution if and only if $u_{i,j}=u_{r,s}$ for
all $1\le i,r, j , s\le t$.
\end{theorem}

We would like now to state a nice conjecture which naturally
arises in this context. Indeed, the results in \cite{R} together
with Theorems \ref{maintheorem}, \ref{theorem2} and \ref{theorem3}
and Corollary \ref{coro1} suggest -and prove in many cases- the
following conjecture

\begin{conjecture}
 Let $I\subset R$ be a graded Cohen-Macaulay ideal of
codimension $c$. Then, we have:
\begin{equation}
 \prod_{1\le j<i}\frac{m_j}{M_i-m_j}\prod_{i<j\le c}\frac{m_j}{M_j-m_i}
\le \beta_i(R/I)\le \frac{1}{(i-1)!·(c-i)!}\prod_{j\ne i}M_j
\end{equation} for $i=1,\cdots ,c$.
\end{conjecture}

 We
will end this paper with some other examples which give support to
the above conjecture.

\begin{theorem} \label{theorem4}
Let $X\subset\PP^n=Proj(S)=Proj(k[x_1, \cdots ,x_{n+1}])$ be a
reduced arithmetically Cohen-Macaulay subscheme of degree
$d>c=codimX$.  Suppose $X$ is a divisor on a variety of minimal
degree.
 Then, we have:
\begin{equation}
 \prod_{1\le j<i}\frac{m_j}{M_i-m_j}\prod_{i<j\le c}\frac{m_j}{M_j-m_i}
\le \beta_i(S/I(X))\le \frac{1}{(i-1)!·(c-i)!}\prod_{j\ne i}M_j
\end{equation} for $i=1,\cdots ,c$.
\end{theorem}
\begin{proof}
First of all recall that $deg(X)=e( S/I(X))$. By \cite{N}, Theorem
1.3, the minimal graded free $S$-resolution of $I(X)$ has the
following shape
$$0\longrightarrow  F_c\longrightarrow F_{c-1}\longrightarrow \cdots
\longrightarrow F_1 \longrightarrow I(X) \longrightarrow 0$$ where
$$F_i=S(-1-i)^{\alpha_i}\oplus S(-t+1-i)^{\beta_i}\oplus
S(-t-i)^{\gamma_i}, \quad 1\le i \le c,$$
$$d=tc+1-p \quad \text{with } 1\le p\le c$$
and $$\alpha_i=i{c\choose i+1}, \quad \text{for } 1\le i\le c,
$$$$\beta_i
= \left \{ \begin{array}{ll}
p{c\choose i-1}-c{c-1\choose i-2} & \hbox{if } 1\le i\le p \\
0 & \mbox{if } p<i\le c,
\end{array} \right .
$$
$$\gamma_i
= \left \{ \begin{array}{ll}
0 & \hbox{if } 1\le i\le p \\
 c{c-1\choose i-1}-p{c\choose i} & \mbox{if } p<i\le c.
\end{array} \right.
$$
Therefore,  we have $$m_i=1+i \quad \text{for } 1\le i\le c,
\text{ and}$$
$$M_i
= \left \{ \begin{array}{ll}
t+i-1 & \hbox{if } 1\le i\le p \\
t+i & \mbox{if } p<i\le c
\end{array} \right.
$$
and the result follows after a straightforward computation taking
into account that $$\beta _{i}(S/I(X))=\alpha _i+\beta _i+\gamma
_i$$ for all all $1\le i \le c$.
\end{proof}

\end{document}